\documentclass[twoside,12pt,a4paper]{article}
\usepackage{amsfonts}
\usepackage{amsmath,amssymb}
\usepackage{amsthm}
\usepackage{graphicx}

\theoremstyle{plain}
\newtheorem{theorem}{Theorem}
\newtheorem{lemma}{Lemma}

\newtheorem*{keywords}{Keywords}
\newtheorem*{amsmos}{AMS(MOS) subject classification}

\title
{
Functional   \hskip 5.8 cm {\small Volume  11     } \\
Differential \hskip 5.4 cm {\small 2004, No 1-2   } \\
Equations    \hskip 6.2 cm  {\small pp.  1--7 \ }  \\
\vskip 1 in
On the simplest system with retarding switching and a
2--point critical set}
\author{D. A. Filimonov \thanks{Department of Mathematics and Mechanics,
Moscow State University, mityafil@yandex.ru}}

\date{}

\begin{document}

\maketitle
\begin{abstract}
The system considered in this paper consists of two equations
$(k=1,2)$ $\dot x(t)=(-1)^{k-1} (0\le t<\infty),\,
k(0)=1,\,x(0)=0,\,x(t)\not\in\{0,1\}(-1\le t<0),$ that change
mutually in every instant $t$ for which $x(t-\tau)\in\{0,1\}$,
where $\tau={\rm const}>0$ is given. In this paper the behavior of
the solutions is characterized for every
$\tau\in(\frac{4}{3},\frac{3}{2})$, i. e. in case not covered in
\cite{ADM}; as it was noted there, this behavior turned out to be
more complex then when $\tau\in(3/2,\infty)$. Thus the behavior of
the solutions of this system with critical set $K=\{0,1\}$ is
characterized for every $\tau>0$.
\end{abstract}

\bigskip

\begin{keywords}
Functional differential equation, Equation of retarded type,
System with retarding switching, Periodic solution
\end{keywords}

\bigskip

\begin{amsmos}
34K15
\end{amsmos}

\bigskip

In \cite{MK} the systems with retarding switching were introduced
as a natural generalization of the "raging systems" of T. Vogel
\cite{TV}. The Vogel's system is defined with two autonomous
systems in ${\bf R^n}$ ("base systems") which replace each other
when the phase point $x(t)$ reaches the given fixed "critical set"
$K$. General systems with retarding switching were defined in
\cite{AM}. Here the change of the base systems occurs in every
instant $t$ when the point $x(t-\tau)$ gets into $K$, where $\tau
={\rm const}>0$ is given, and the values of the solution are given
on some time interval of the length $\tau$ as the initial
condition, as well as the number of the initial base system.

In \cite{ADM} the considered system has the form
\begin{equation}
\label{1}
\dot x(t)=(-1)^{k-1}\qquad(k=1,2;\quad0\le t<\infty),
\end{equation}
whereas the critical set is $K=\{0,1\}$. We set the continuous
initial function $x=\varphi(t)$, for $-1\le t\le0$, where $
\varphi(0)=0,\,\varphi(t)\not\in\{0,1\}(-1\le t<0)$, and start
with the first base equation $(k=1)$.

In \cite{ADM} the behavior of the solution of the system (\ref{1})
is fully characterized for
$\tau\in[0,\frac{4}{3}]\cup[\frac{3}{2},\infty)$. When
$\tau\in[\frac{4}{3},\frac{3}{2}]$ there are some calculations for
$\tau\in[\frac{4}{3},\frac{31}{21}]$.

But one critical value - $\tau=\frac{63}{43}$ - in interval
$(\frac{16}{11},\frac{31}{21})$ was lost. That's why coordinates
$x$ for turning points mentioned in \cite{ADM} are correct only
for $\tau\in(\frac{16}{11}, \frac{63}{43})$.

If $\tau=\frac{63}{43}$ after switching in points $\frac{63}{43},
\frac{20}{43}, \frac{60}{43}, \frac{14}{43}, \frac{48}{43},
-\frac{10}{43}, 0, -1, -\frac{23}{43}$ the solution goes to
$-\infty$;

if $\tau\in(\frac{63}{43},\frac{31}{21})$ it is periodic with
sequential turning points
$\tau,\tau-1,3\tau-3,5\tau-7,11\tau-15,21\tau-31,43\tau-63,
43\tau-64,\tau-2,-85\tau+124$.

Probably this oversight hindered the author of \cite{ADM} from
finding the general rule of behavior of the solutions of the
system (\ref{1}) when $\tau\in[4/3,3/2)$.

In this paper the behavior of the solutions of the system
(\ref{1}) is characterized for every
$\tau\in(\frac{4}{3},\frac{3}{2})$, i. e. in case not covered in
\cite{ADM}; as it was noted there, this behavior turned out to be
more complex then when $\tau\in(3/2,\infty)$. Thus the behavior of
the solution of this system with critical set $K=\{0,1\}$ is
characterized for every $\tau>0$.

Let us introduce the designations
$\tau_k:=3\cdot4^{k}/(2\cdot4^k+1)\enspace(k\in\mathbb N)$;
$\theta_k:=3\cdot(4^{k+1}-1)/(2\cdot4^{k+1}+1)\enspace(k\in\mathbb
N)$ and
$\zeta_k:=3\cdot(2\cdot4^k-1)/(4^{k+1}-1)\enspace(k\in\mathbb N)$.
This are the increasing  sequences, with $\tau_1=\frac{4}{3},
\tau_k\to3/2$ as $k\to\infty$; $\theta_1=\frac{15}{11},
\theta_k\to3/2$ as $k\to\infty$ and $\zeta_1=\frac{7}{5},
\zeta_k\to3/2$ as $k\to\infty$. More than
$\tau_k<\theta_k<\zeta_k<\tau_{k+1}\,(\forall k\in\mathbb N)$.

\begin{theorem}
\label{T}

For all $k\in\mathbb N$:

the solution of the problem is periodic and has $4k+2$ switchings
on the least period, if $\tau=\tau_k$;

the solution of the problem is periodic and has $2k+4$ switchings
on the least period, if $\tau_{k}<\tau<\theta_k$;

the solution of the problem goes to $-\infty$ after $2k+5$
switchings, if $\tau=\theta_k$;

the solution of the problem is periodic and has $2k+6$ switchings
on the least period, if $\theta_{k}<\tau<\zeta_k$;

the solution of the problem goes to $-\infty$ after $4k+5$
switchings, if $\tau=\zeta_k$;

the solution of the problem is periodic and has $2k+4$ switchings
on the least period, if $\zeta_{k}<\tau<\tau_{k+1}$.

\end{theorem}

\proof
Let consider, that $\tau\in[4/3,3/2)$ is given. Let
$(a_1,b_1)(=(0,0)),\\(a_2,b_2),\ldots,\,((\alpha_1,\beta_1),(
\alpha_2,\beta_2),\ldots)$ be the sequence of the values of $(x,
t)$ in the instants of the hit of the solution on the critical set
(in instants of the switching, respectively).

If $\tau=\tau_1=\frac{4}{3}$ the solution is periodic with
sequential turning points
$\frac{4}{3},\frac{1}{3},1,-\frac{1}{3},\\ \frac{2 }{3},0$ and it
corresponds with the formulation of the Theorem \ref{T}.

Let us prove that Theorem \ref{T} is valid for all
$\tau\in(4/3,3/2)$. Denote by $J$ the maximal natural number for
which $ \alpha_j>1$ at all odd $j<J$ and $\alpha_j<1$ at all even
$j<J$; $J=\infty$ if $\alpha_j$ satisfy these conditions for all
natural $j$. It is obvious  that $J\ge4$.

\begin{lemma}
\label{L}
Let $\tau\in(4/3,3/2)$, $2k+1<J$; then $\alpha_{2k}>0$.
\end{lemma}

\proof
Let proof this statement by contradiction. Let assume that
there exists $k\in\mathbb N$, so that $2k+1<J$ but
$\alpha_{2k}\le0$; from all such $k$ choose the smallest. From the
definition of $J$ and from our assumption follows, that for all
$n\in\mathbb{N}\enspace (1\le n\le2k-2)$ in every segment with
ends $\alpha_n, \alpha_{n+1}$ there is only one point from our
critical set : $\{1\}$. Taking into account formulas
$a_2=1,\enspace \beta_1<b_3<\beta_2$, we conclude, that
$a_n=1\,(1\le n\le2k)$ and $\beta_{n-2}<b_n<\beta_{n-1}\,(3\le
n\le2k+1)$. From $\alpha_{2k-1}>0$ and $\alpha_{2k}\le0$ we
obtain, that $a_{2k+1}=1,\enspace a_{2k+2}=0$. According to the
definition of solution we have $\beta_n-\beta_n=\tau (\forall
n\in\mathbb N)$; in particular, $\beta_{2k+1}-b_{2k+1}=\tau<3/2$.
But $\beta_{2k+1}-b_{2k+1}=
(\beta_{2k+1}-\beta_{2k})_(\beta_{2k}-b_{2k+1}=|\alpha_{2k+1}-
\alpha_{2k}|+|\alpha_{2k}-a_{2k}|>2$. This contradiction shows us
falseness of our assumption. Lemma \ref{L} is proved. 
\endproof

Let $j\in[2,J-1]$. It follows from Lemma \ref{L}, that
$\alpha_{j}>0$ for all $j\in[2,J-2]$. Then it follows from the
definition of the solution that
\begin{equation}
\label{2}
b_{j+1}-\beta_{j-1}=\beta_{j-1}-b_j.
\end{equation}
But $\beta_j=b_j+\tau$ for all values of $j$ such that $b_j$ is
defined - in particular, for $j=J<\infty$, as $b_J<\beta_{J-1}$.
Therefore from (\ref{2})  we obtain the recurrence relation
\begin{equation}
\label{3}
\beta_{j+1}=-\beta_j+2\beta_{j-1}+2\tau,\quad2\le j\le
J-1.
\end{equation}
Moreover the definition of the solution directly implies the
expression for $\alpha$ in terms of $\beta$:
\begin{equation}
\label{4}
\alpha_j=1+(-1)^j[\tau-2(\beta_j-\beta_{j-1})],\quad2\le
j\le J.
\end{equation}

The formula
$$\beta_j=\frac{6j+1-(-2)^j}{9}\,\tau-\frac{(-2)^{j-1}-1}{3}\,,
\quad1\le j\le J,$$ follows from the relation (\ref{3}) and the
initial data $\beta_1=\tau ,\,\beta_2=\tau+1$. This formula
together with (\ref{4}) and the equality $\alpha_1=\tau$ yields
that
\begin{equation}
\label{5} \alpha_j=\frac{2^j-(-1)^j}{3}\,\tau-2^{j-1}+1,\quad1\le
j\le J.
\end{equation}
One can see that both sequences $\{\alpha_j\}(j\in[1,J]$ even) and
$(j\in[1,J]$ odd) are decreasing. Hence $J$ is the smallest even
number for which the value $\alpha_j$ calculated according to the
formula (\ref{5}), becomes $\le1$. The values $\tau_k,k\in\mathbb
N$ for which $\alpha_{2k+1}=1$, i.e.
$$\frac{2^{2k+1}+1}{3}\tau_k-2^{2k}+1=1.$$
are the critical ones. We find from here that
$\tau_k=3\cdot4^{k}/(2\cdot4^k+1)$. Now let's consider possible
cases.

If $\tau=\tau_k$, then $\alpha_{2k+1}=1$, that's why because in
$[\alpha_{2k},\alpha_{2k+1}]$ there are no other points from
critical set, the solution continues with the same time intervals
but symmetrically relatively $\{x=1/2\}$, i. e. the solution is
periodical and has $2k+1+2k+1=4k+2$ switchings on the least period
(Fig. 1)\footnote{All figures were rendered for $k=3$}.

\unitlength=0.7cm

\begin{picture}(18,6)(0.3, -3)
\put(0,0){\vector(1,0){18}} \put(17.7, -0.5){$t$}
\put(0,0){\vector(0,1){2}} \put(-0.5,1.7){$x$} \put(-0.5,0.8){$1$}
\put(-0.5,-0.2){$0$} \put(0,1){\line(1,0){18}}
\put(0.000000000000000000,0.000000000000000000){\line(1,1){1.488372093023255810}}
\put(1.488372093023255810,1.488372093023255810){\line(1,-1){1.000000000000000000}}
\put(2.488372093023255810,0.488372093023255814){\line(1,1){0.976744186046511628}}
\put(3.465116279069767440,1.465116279069767440){\line(1,-1){1.023255813953488370}}
\put(4.488372093023255810,0.441860465116279070){\line(1,1){0.930232558139534883}}
\put(5.418604651162790700,1.372093023255813950){\line(1,-1){1.116279069767441860}}
\put(6.534883720930232560,0.255813953488372092){\line(1,1){0.744186046511627906}}
\put(6.2,-0.5){$\alpha_{2k}$}
\put(7.279069767441860460,0.999999999999999998){\line(1,-1){1.488372093023255820}}
\put(6.9,1.3){$\alpha_{2k+1}$}
\put(8.767441860465116280,-0.488372093023255818){\line(1,1){0.999999999999999996}}
\put(9.767441860465116280,0.511627906976744179){\line(1,-1){0.976744186046511635}}
\put(10.744186046511627900,-0.465116279069767456){\line(1,1){1.023255813953488360}}
\put(11.767441860465116300,0.558139534883720901){\line(1,-1){0.930232558139534912}}
\put(12.697674418604651200,-0.372093023255814011){\line(1,1){1.116279069767441800}}
\put(13.813953488372093000,0.744186046511627792){\line(1,-1){0.744186046511628022}}
\put(14.558139534883721000,-0.000000000000000229){\line(1,1){1.488372093023255590}}
\put(14,-0.5){$\alpha_{4k+2}$} \put(8,-2.5){Fig. 1.}
\end{picture}

Let then $\tau_k<\tau<\tau_{k+1}$. Therefore $\alpha_{2k+1}>1$ and
$\alpha_{2k+3}<1$. The further behavior of the solution depends on
that, if $\alpha_{2k+3}$ is greater than zero or no. So we have
another set of critical points, for which
$$\frac{2^{2k+3}+1}{3}\theta_k-2^{2k+2}+1=0,$$ wherefrom we
receive $\theta_k=3\cdot(4^{k+1}-1)/(2\cdot4^{k+1}+1)$.

If $\tau_k<\tau<\theta_k$, i. e. $\alpha_{2k+3}<0$, then after
switchings in points $\alpha_{2k+2}$ and $\alpha_{2k+3}$ there is
one another switching in point $\alpha_{2k+4}<0$ and the solution
gets to the beginning of it's period after $2k+4$ switchings (Fig.
2).

\begin{picture}(18,6)(0.3, -3)
\put(0,0){\vector(1,0){18}} \put(17.7, -0.5){$t$}
\put(0,0){\vector(0,1){2}} \put(-0.5,1.7){$x$} \put(-0.5,0.8){$1$}
\put(-0.5,-0.2){$0$} \put(0,1){\line(1,0){18}}
\put(6.2,-0.5){$\alpha_{2k}$} \put(6.9,1.3){$\alpha_{2k+1}$}
\put(8.4,-0.9){$\alpha_{2k+2}$} \put(8.9,0.2){$\alpha_{2k+3}$}
\put(9.9,-1.4){$\alpha_{2k+4}$}
\put(0.000000000000000000,0.000000000000000000){\line(1,1){1.490372093023255810}}
\put(1.490372093023255810,1.490372093023255810){\line(1,-1){1.000000000000000000}}
\put(2.490372093023255810,0.490372093023255814){\line(1,1){0.980744186046511628}}
\put(3.471116279069767440,1.471116279069767440){\line(1,-1){1.019255813953488370}}
\put(4.490372093023255810,0.451860465116279070){\line(1,1){0.942232558139534884}}
\put(5.432604651162790700,1.394093023255813950){\line(1,-1){1.096279069767441860}}
\put(6.528883720930232560,0.297813953488372093){\line(1,1){0.788186046511627907}}
\put(7.317069767441860470,1.086000000000000000){\line(1,-1){1.704372093023255810}}
\put(9.021441860465116280,-0.618372093023255814){\line(1,1){0.542000000000000000}}
\put(9.563441860465116280,-0.076372093023255814){\line(1,-1){1.000000000000000000}}
\put(10.563441860465116280,-1.076372093023255810){\line(1,1){1.146372093023255810}}
\put(11.709813953488372100,0.070000000000000000){\line(1,1){1.490372093023255810}}
\put(8,-2.5){Fig. 2.}
\end{picture}

If $\tau=\theta_k$, i. e. $\alpha_{2k+3}=0$, then we can verify by
direct calculating, that $\alpha_{2k+4}=-1$. Therefore
$\alpha_{2k+5}<0$ after which nothing does not prevent the
solution to go to $-\infty$. It has performed $2k+5$ switchings
(Fig. 3).

\begin{picture}(18,6)(0.3, -3)
\put(0,0){\vector(1,0){18}} \put(17.7, -0.5){$t$}
\put(0,0){\vector(0,1){2}} \put(-0.5,1.7){$x$} \put(-0.5,0.8){$1$}
\put(-0.5,-0.2){$0$} \put(0,1){\line(1,0){18}}
\put(6.2,-0.5){$\alpha_{2k}$} \put(6.9,1.3){$\alpha_{2k+1}$}
\put(8.2,-1){$\alpha_{2k+2}$} \put(8.9,0.4){$\alpha_{2k+3}$}
\put(9.9,-1.5){$\alpha_{2k+4}$} \put(10.5,-0.3){$\alpha_{2k+5}$}
\put(0.000000000000000000,0.000000000000000000){\line(1,1){1.491228070175438600}}
\put(1.491228070175438600,1.491228070175438600){\line(1,-1){1.000000000000000000}}
\put(2.491228070175438600,0.491228070175438597){\line(1,1){0.982456140350877193}}
\put(3.473684210526315790,1.473684210526315790){\line(1,-1){1.017543859649122810}}
\put(4.491228070175438600,0.456140350877192983){\line(1,1){0.947368421052631579}}
\put(5.438596491228070180,1.403508771929824560){\line(1,-1){1.087719298245614030}}
\put(6.526315789473684210,0.315789473684210527){\line(1,1){0.807017543859649124}}
\put(7.333333333333333330,1.122807017543859650){\line(1,-1){1.668421052631578950}}
\put(9.001754385964912280,-0.545614035087719294){\line(1,1){0.545614035087719302}}
\put(9.547368421052631580,0.000000000000000008){\line(1,-1){1.000000000000000000}}
\put(10.547368421052631580,-0.999999999999999992){\line(1,1){0.591228070175438588}}
\put(11.138596491228070200,-0.408771929824561404){\line(1,-1){1}}
\put(8,-2.5){Fig. 3.}
\end{picture}

If $\theta_k<\tau<\tau_{k+1}$, the result depends on the sign of
$\alpha_{2k+2}$. From this we receive the last set of critical
points $$\frac{2^{2k+2}-1}{3}\zeta_k-2^{2k+1}+1=0,$$ or
$\zeta_k=3\cdot(2\cdot4^k-1)/(4^{k+1}-1)$.

If $\theta_k<\tau<\zeta_k$, then we can verify by direct
calculating, that $\alpha_{2k+5}=\tau-2<0$. So the solution after
switching in points $\alpha_{2k+4}$, $\alpha_{2k+5}$ and
$\alpha_{2k+6}$ gets to the beginning of it's period after $2k+6$
switchings (Fig. 4).

\begin{picture}(18,6)(0.3, -3)
\put(0,0){\vector(1,0){18}} \put(17.7, -0.5){$t$}
\put(0,0){\vector(0,1){2}} \put(-0.5,1.7){$x$} \put(-0.5,0.8){$1$}
\put(-0.5,-0.2){$0$} \put(0,1){\line(1,0){18}}
\put(6.2,-0.5){$\alpha_{2k}$} \put(6.9,1.3){$\alpha_{2k+1}$}
\put(8.1,-0.6){$\alpha_{2k+2}$} \put(8.7,0.5){$\alpha_{2k+3}$}
\put(9.5,-1.0){$\alpha_{2k+4}$} \put(10.12,0.2){$\alpha_{2k+5}$}
\put(10.9,-1.0){$\alpha_{2k+6}$}
\put(0.000000000000000000,0.000000000000000000){\line(1,1){1.492228070175438600}}
\put(1.492228070175438600,1.492228070175438600){\line(1,-1){1.000000000000000000}}
\put(2.492228070175438600,0.492228070175438597){\line(1,1){0.984456140350877193}}
\put(3.476684210526315790,1.476684210526315790){\line(1,-1){1.015543859649122810}}
\put(4.492228070175438600,0.461140350877192983){\line(1,1){0.953368421052631579}}
\put(5.445596491228070180,1.414508771929824560){\line(1,-1){1.077719298245614030}}
\put(6.523315789473684210,0.336789473684210528){\line(1,1){0.829017543859649124}}
\put(7.352333333333333340,1.165807017543859650){\line(1,-1){1.436421052631578940}}
\put(8.788754385964912280,-0.270614035087719292){\line(1,1){0.551614035087719305}}
\put(9.340368421052631580,0.281000000000000014){\line(1,-1){1.000000000000000000}}
\put(10.340368421052631600,-0.718999999999999986){\line(1,1){0.541228070175438583}}
\put(10.881596491228070200,-0.177771929824561403){\line(1,-1){0.542000000000000028}}
\put(11.423596491228070200,-0.719771929824561431){\line(1,1){0.849771929824561431}}
\put(12.273368421052631600,0.130000000000000000){\line(1,1){1.492228070175438600}}
\put(8,-2.5){Fig. 4.}
\end{picture}

If $\tau=\zeta_k$. Then we can observe, that
$\alpha_{2k+3}=\tau-1$. It means, that the solution further behave
as if $t\in[\beta_1,\beta_{2k+3}]$, and after that, having in all
$2k+2+2k+2+1=4k+5$ switchings, goes to $-\infty$ (Fig. 5).

\begin{picture}(18,6)(0.3, -3)
\put(0,0){\vector(1,0){18}} \put(17.7, -0.5){$t$}
\put(0,0){\vector(0,1){2}} \put(-0.5,1.7){$x$} \put(-0.5,0.8){$1$}
\put(-0.5,-0.2){$0$} \put(0,1){\line(1,0){18}}
\put(6.2,-0.5){$\alpha_{2k}$} \put(6.9,1.5){$\alpha_{2k+1}$}
\put(7.9,-0.5){$\alpha_{2k+2}$} \put(8.5,0.7){$\alpha_{2k+3}$}
\put(16.3,0.2){$\alpha_{4k+5}$}
\put(0.000000000000000000,0.000000000000000000){\line(1,1){1.494117647058823530}}
\put(1.494117647058823530,1.494117647058823530){\line(1,-1){1.000000000000000000}}
\put(2.494117647058823530,0.494117647058823529){\line(1,1){0.988235294117647059}}
\put(3.482352941176470590,1.482352941176470590){\line(1,-1){1.011764705882352940}}
\put(4.494117647058823530,0.470588235294117646){\line(1,1){0.964705882352941176}}
\put(5.458823529411764700,1.435294117647058820){\line(1,-1){1.058823529411764710}}
\put(6.517647058823529410,0.376470588235294115){\line(1,1){0.870588235294117644}}
\put(7.388235294117647060,1.247058823529411760){\line(1,-1){1.247058823529411770}}
\put(8.635294117647058830,-0.000000000000000011){\line(1,1){0.594117647058823519}}
\put(9.229411764705882350,0.594117647058823508){\line(1,-1){1.000000000000000000}}
\put(10.229411764705882400,-0.405882352941176492){\line(1,1){0.988235294117647038}}
\put(11.217647058823529400,0.582352941176470547){\line(1,-1){1.011764705882352980}}
\put(12.229411764705882400,-0.429411764705882436){\line(1,1){0.964705882352941094}}
\put(13.194117647058823500,0.535294117647058658){\line(1,-1){1.058823529411764870}}
\put(14.252941176470588300,-0.523529411764706214){\line(1,1){0.870588235294117316}}
\put(15.123529411764705600,0.347058823529411102){\line(1,-1){1.247058823529412430}}
\put(16.370588235294118100,-0.900000000000001330){\line(1,1){0.594117647058822204}}
\put(16.964705882352940300,-0.305882352941179121){\line(1,-1){1}}
\put(8,-2.5){Fig. 5.}
\end{picture}

And finally if $\zeta_k<\tau<\tau_{k+1}$, then because in this
case $\alpha_{2k+2}>0$ and $\alpha_{2k+3}<1$, there is only one
another switching, after which the solution gets to the beginning
of it's period after $2k+4$ switchings (Fig. 6). This ends proof
of Theorem 1.

\begin{picture}(18,6)(0.3, -3)
\put(0,0){\vector(1,0){18}} \put(17.7, -0.5){$t$}
\put(0,0){\vector(0,1){2}} \put(-0.5,1.7){$x$} \put(-0.5,0.8){$1$}
\put(-0.5,-0.2){$0$} \put(0,1){\line(1,0){18}}
\put(6.2,-0.5){$\alpha_{2k}$} \put(6.9,1.6){$\alpha_{2k+1}$}
\put(7.9,-0.5){$\alpha_{2k+2}$} \put(8.6,1.2){$\alpha_{2k+3}$}
\put(11,-2){$\alpha_{2k+4}$}
\put(0.000000000000000000,0.000000000000000000){\line(1,1){1.496117647058823530}}
\put(1.496117647058823530,1.496117647058823530){\line(1,-1){1.000000000000000000}}
\put(2.496117647058823530,0.496117647058823529){\line(1,1){0.992235294117647059}}
\put(3.488352941176470590,1.488352941176470590){\line(1,-1){1.007764705882352940}}
\put(4.496117647058823530,0.480588235294117647){\line(1,1){0.976705882352941176}}
\put(5.472823529411764710,1.457294117647058820){\line(1,-1){1.038823529411764710}}
\put(6.511647058823529410,0.418470588235294116){\line(1,1){0.914588235294117645}}
\put(7.426235294117647060,1.333058823529411760){\line(1,-1){1.163058823529411770}}
\put(8.589294117647058830,0.169999999999999993){\line(1,1){0.666117647058823523}}
\put(9.255411764705882350,0.836117647058823516){\line(1,-1){0.836117647058823516}}
\put(10.091529411764705900,0.000000000000000000){\line(1,-1){1.496117647058823530}}
\put(11.587647058823529400,-1.496117647058823530){\line(1,1){1.496117647058823530}}
\put(13.083764705882352900,0.000000000000000000){\line(1,1){1.496117647058823530}}
\put(8,-2.5){Fig. 6.}
\end{picture}

\endproof

{\bf Acknowledgement.} I thank Prof. A.D.Myshkis for the setting
of a problem and for useful discussion.

\end{document}